\documentclass[10pt]{amsart}
\usepackage{amssymb, amsthm}\usepackage{amsmath}\usepackage{epsfig}
\usepackage{fancyhdr}

\theoremstyle{plain}
\newtheorem{Lem}{Lemma}[section]
\newtheorem{Not}[Lem]{Notation}
\newtheorem{Prop}[Lem]{Proposition}
\newtheorem{Cor}[Lem]{Corollary}

\newtheorem{The}[Lem]{Theorem}

\theoremstyle{definition}

\newtheorem{definition}[Lem]{Definition}
\newtheorem{Rem}[Lem]{Remark}
\newtheorem{Exe}[Lem]{Example}
\def\sn{W}
\def\ZZe{\sn(e)}
\def\ZZf{\sn(f)}
\def\ZZse{\sn_\star(e)}
\def\ZZie{\sn^\star(e)}
\def\De{D(e)}
\def\Ge{G(e)}
\def\Gue{G^\uparrow(e)}

\def\Df{D(f)}

\def\Duf{D^\uparrow(f)}

\def\Dh{D(h)}
\def\Gh{G(h)}

\def\Guei{G^\uparrow(e_i)}
\def\Duei{D^\uparrow(e_i)}

\def\Gueun{G^\uparrow(e_1)}
\def\Duej{D^\uparrow(e_j)}
\def\Guel{G^\uparrow(e_\ell)}
\def\Duel{D^\uparrow(e_\ell)}
\def\Guell{G^\uparrow(e_{\ell-1})}
\def\Duell{D^\uparrow(e_{\ell-1})}
\def\Gufl{G^\uparrow(f_\ell)}
\def\Dufl{D^\uparrow(f_\ell)}
\def\Dufun{D^\uparrow(f_1)}

\def\Dse{D_\star(e)}
\def\Gse{G_\star(e)}
\def\Dsh{D_\star(h)}
\def\Gsh{G_\star(h)}
\def\renn{R(M)}
\def\crl{\Lambda}
\def\crlo{\Lambda_{\scriptscriptstyle \circ}}
\begin{document}
\date{\today}
\author{Godelle Eddy}
\title{A note on Renner monoids}
\maketitle
\begin{abstract}
We extend the result obtained in~\cite{God} to every Renner monoid: we provide Renner monoids with a monoid presentation and we introduce a length function which extends the Coxeter length function and which behaves nicely. 
\end{abstract}
%--------------------------------------------------
%--------------------------------------------------
\section*{Introduction}
%--------------------------------------------------

The notion of a \emph{Weyl group} is crucial in Linear Algebraic Group Theory~\cite{Hum}. The seminal example occurs when one considers the \emph{algebraic group}~$GL_n(K)$. In that case, the associated Weyl group is isomorphic to the group of monomial matrices, that is to the permutation group~$S_n$. Weyl groups are special examples of \emph{finite Coxeter groups}. Hence, they possess a group presentation of a particular type, and an associated length function. It turns out that this presentation and this length function are deeply related with the geometry of the associated algebraic group. Linear Algebraic Monoid Theory, mainly developed by Putcha, Renner and Solomon, has deep connections with  Algebraic Group Theory. In particular, the \emph{Renner monoid}~\cite{Ren} plays the role that the Weyl group does in Linear Algebraic Group Theory. As far as I know, in the case of Renner monoids, there is no known theory that plays the role of Coxeter Group Theory. Therefore it is natural to look for such a theory, and, therefore, to address the question of monoid presentations for Renner monoids. In~\cite{God}, we considered the particular case of the \emph{rook monoid} defined by Solomon~\cite{Sol4}. We obtained a presentation of this group and introduced a length function, that is nicely related to the Hecke Algebra of the rook monoid. Our objective here is to consider the general case. We obtain a presentation of every Renner monoid and introduce a length function. In the case of the rook monoid, we recover the results obtained in~\cite{God}. Our length function is not the classical length function on Renner monoids~\cite{Ren}. We remark that the former shares with the latter several nice geometrical and combinatorial properties. We believe that this new length function is nicely related to the Hecke algebra, as in the special case of the rook monoid.

Let us postpone to the next sections some definitions and notations, and state here our main results. 
Consider the Renner monoid~$\renn$ of a linear algebraic monoid~$M$. Denote by $W$ the unit group of $\renn$ and consider its associated Coxeter system~$(W,S)$. Denote by $\crl$ a \emph{cross section lattice} of the monoid~$E(R(M))$ of idempotent elements of $R(M)$, and by $\crlo$ the set of elements of~$\crl$ that are distinct from the unity. Denote finally by $\lambda$ the associated \emph{type map} of $\renn$; roughly speaking, this is a map that describes the action of $W$ on $E(R(M))$.

\begin{The} The Renner monoid~$\renn$ admits the monoid presentation whose generating set is~$S\cup \crlo$ and whose defining relations are:\\  
\begin{tabular}{lll}
(COX1)&$s^2 = 1$,&$s\in S$;\\
(COX2)&$|s,t\rangle^m = |t,s\rangle^m$,&$(\{s,t\},m)\in \mathcal{E}(\Gamma)$;\\
(TYM1)&$se = es$,&  $e\in \crlo$, $s\in \lambda^\star(e)$;\\
(TYM2)&$se = es = e$,& $e\in \crlo$, $s\in \lambda_\star(e)$;\\
(TYM)&$e\underline{w}f = e\wedge_wf$,& $e,f\in \crlo$, $w\in \Gue\cap \Duf$. 
\end{tabular} 
\label{thepreserennmonidintro}
\end{The}

We define the length~$\ell$ on~$\renn$ in the following way: if~$s$ lies in~$S$, we set~$\ell(s) = 1$; if $e$ lies in $\crl$ we set~$\ell(e) = 0$. Then, we extend $\ell$ by additivity to the free monoid of words on~$S\cup \crlo$. If $w$ lies in~$\renn$, its length~$\ell(w)$  is the minimal length of its word representatives on~$S\cup \crlo$. In Section~2, we investigate the properties of this length function. In particular we prove that it is nicely related to the classical \emph{normal form} defined on~$\renn$, and we also prove  
\begin{Prop}\label{propintro} Let $T$ be a maximal torus of the unit group of~$M$. Fix a Borel subgroup~$B$ that contains $T$. Let~$w$ lie in~$\renn$ and~$s$ lie in~$S$. Then, $$B s B w B = \left\{\begin{array}{ll}BwB&\textrm{if } \ell(sw) = \ell(w);\\BswB&\textrm{if } \ell(sw) = \ell(w)+1;\\BswB\cup BwB&\textrm{if } \ell(sw) = \ell(w)-1.\end{array}\right.$$\end{Prop}

This article is organized as it follows. In section~$1$, we first recall the backgrounds on Algebraic Monoid Theory and on Coxeter Theory. Then, we prove Theorem~\ref{thepreserennmonidintro}. In Section~2, we  consider several examples of Renner monoids and deduce explicit presentations from Theorem~\ref{thepreserennmonidintro}. In Section~3, we focus on the length function.
%--------------------------------------------------
%
%--------------------------------------------------
%--------------------------------------------------
\section{Presentation for Renner monoids}
%--------------------------------------------------
%--------------------------------------------------
%
Our objective in the present section is to associate a monoid presentation to every Renner monoid. The statement of our result and its proof require some properties of algebraic monoid theory and of Coxeter group theory. In Section~\ref{sousectionrappelamt}, we introduce Renner monoids and we state the results we need about algebraic monoids. In Section~\ref{sousectionrappelcgt} we recall the definition of Coxeter groups and some of its well-known properties. Using the two preliminary sections, we can prove~Theorem~\ref{thepreserennmonidintro} in Section~\ref{sousectionpresrenmnoid}. This provides a monoid presentation to every Renner monoid.     

We fix an algebraically closed field~$\mathbb{K}$. We denote by~$M_n$ the set of all~$n\times n$ matrices over~$\mathbb{K}$, and by~$GL_n$ the set of all invertible matrices in~$M_n$. We refer to~\cite{Put,Ren,Sol1} for the general theory and proofs involving linear algebraic monoids and Renner monoids; we refer to~\cite{Hum} for an introduction to linear algebraic groups. If $X$ is a subset of $M_n$, we denote by $\overline{X}$  its closure for the Zariski topology.
 
%--------------------------------------------------
\subsection{Algebraic Monoid Theory}
%--------------------------------------------------
\label{sousectionrappelamt}
 We introduce here the basic definitions and notation on Algebraic Monoid Theory that we shall need in the sequel.  
%--------------------------------------------------
\subsubsection{Regular monoids and reducible groups}
%--------------------------------------------------
\begin{definition}[Algebraic monoid] An \emph{algebraic monoid} is a submonoid of~$M_n$, for some positive integer~$n$, that is closed for the Zariski topology. An algebraic monoid is \emph{irreducible} if it is irreducible as a variety.\end{definition}

It is very easy to construct algebraic monoids. Indeed, the Zariski closure~$M = \overline{G}$ of any submonoid~$G$ of $M_n$ is an algebraic monoid. The main example occurs when for $G$ one considers an algebraic subgroup of~$GL_n$. It turns out that in this case, the group~$G$ is the unit group of~$M$. Conversely, if $M$ is an algebraic monoid, then its unit group~$G(M)$ is an algebraic group. The monoid~$M_n$ is the seminal example of an algebraic monoid, and its unit group~$GL_n$ is the seminal example of an algebraic group.

One of the main differences between an algebraic group and an algebraic monoid is that the latter have idempotent elements. In the sequel we denote by~$E(M)$ the set of idempotent elements of a monoid~$M$.  We recall that $M$ is \emph{regular} if $M = E(M)G(M) = G(M)E(M)$, and that $M$ \emph{has a zero element} if there exists an element~$0$ such that~$0\times m = m\times 0 = 0$ for every~$m$ in~$M$. The next result, which is the starting point of the theory, was obtained independently by Putcha and Renner in 1982. 
\begin{The} Let $M$ be an irreducible algebraic monoid with a zero element. Then $M$ is regular if and only if $G(M)$ is reductive.  
\end{The}

The order~$\leq$ on~$E(M)$ defined by~$e\leq f$ if~$ef =fe = e$, provides a natural connection between the Borel subgroups of $G(M)$ and the idempotent elements of~$M$ : 
\begin{The} \label{theoemoppborel}Let $M$ be a regular irreducible algebraic monoid with a zero element. Let~$\Gamma = (e_1,\ldots, e_k)$ be a maximal increasing sequence of distinct elements of~$E(M)$.\\ (i) The centralizer~$Z_{G(M)}(\Gamma)$ of~$\Gamma$ in~$G(M)$ is a maximal torus of~the reducive group~$G(M)$.\\(ii) Set $$B^+(\Gamma) = \{b\in G(M)\mid \forall e\in \Gamma, be = ebe\},$$  $$B^-(\Gamma) = \{b\in G(M)\mid \forall e\in \Gamma, eb = ebe\}.$$ Then, $B^-(\Gamma)$ and $B^+(\Gamma)$ are two opposed Borel subgroups that contain~$Z_{G(M)}(\Gamma)$.
\end{The}  

%----------------------------------------
\subsubsection{Renner monoid}
%----------------------------------------
\begin{definition}[Renner monoid] Let $M$ be a regular irreducible algebraic monoid with a zero element. If $T$ is a Borel subgroup of $G(M)$, then we denote its normalizer by $N_{G(M)}(T)$. The \emph{Renner monoid}~$\renn$ of~$M$ is the monoid~$\overline{N_{G(M)}(T)}/T$. 
\end{definition}
It is clear that $\renn$ does not depend on the choice of the maximal torus of $G(M)$.
\begin{Exe} Consider $M = M_n(K)$, and choose the maximal torus~$\mathbb{T}$ of diagonal matrices. The Renner monoid is isomorphic to the monoid of matrices with at most one nonzero entry, that is equal to~$1$, in each row and each column. This monoid is called the rook monoid~$R_n$~\cite{Sol2}. Its unit group is the group of monomial matrices, which is isomorphic to the symmetric group~$S_n$.\end{Exe}

It is almost immediate from the definition that we have
\begin{Prop} Let $M$ be a regular irreducible algebraic monoid with a zero element, and fix a maximal torus~$T$ of $G(M)$. The Renner monoid~$R(M)$ is a finite factorisable inverse monoid. In particular, the set $E(R(M))$ is a commutative monoid and a lattice for the partial order $\leq$ defined by $e\leq f$ when $ef = e$. Furthermore, there is a canonical order preserving isomorphism of monoids between $E(R(M))$  and $E(\overline{T})$. \label{propRgenparEetW}
\end{Prop}

%--------------------------------------------------
\subsection{Coxeter Group Theory}
%-------------------------------------------------- 
\label{sousectionrappelcgt}
Here we recall some well-known facts about Coxeter groups. We refer to~\cite{Bou} for general theory and proof.
\begin{definition} Let $\Gamma$ be a finite simple labelled graph whose labels are positive integers greatest or equal than~$3$. We denote by $S$ the vertex set of $\Gamma$. We denote by $\mathcal{E}(\Gamma)$ the set of pairs $(\{s,t\},m)$ such that either $\{s,t\}$ is an edge of $\Gamma$ labelled by~$m$ or $\{s,t\}$ is not an edge of $\Gamma$ and $m=2$. When $(\{s,t\},m)$ belongs to~$\mathcal{E}(\Gamma)$, we denote by $|s,t\rangle^m$ the word~$sts\cdots$ of length~$m$. The \emph{Coxeter group}~$W(\Gamma)$ associated with~$\Gamma$ is defined by the following group presentation $$\left\langle S \left| \begin{array}{ll}s^2 = 1&s\in S\\ |s,t\rangle^m = |t,s\rangle^m &(\{s,t\},m)\in \mathcal{E}(\Gamma) \end{array}\right.\right\rangle$$
In this case, one says that $(W,S)$ is a \emph{Coxeter system}.\end{definition}

\begin{Prop}  Let $M$ be a regular irreducible algebraic monoid with a zero element, and denote by $G$ its unit group. Fix a maximal torus~$T$ included in~$B$. Then\\(i) the \emph{Weyl group}~$W = N_G(T)/T$ of $G$ is a finite Coxeter group.\\(ii) The unit group of~$R(M)$ is the Weyl group~$W$. \label{propgather2}
\end{Prop}
\begin{Rem}\label{remgenrset}
Gathering the results of Propositions~\ref{propRgenparEetW} and~\ref{propgather2} we get $$R(M) = E(\overline{T})\cdot W = W\cdot E(\overline{T}).$$ 
\end{Rem}

\begin{definition} Let~$(W,S)$ be a \emph{Coxeter system}. Let $w$ belong to $\sn$. The \emph{length}~$\ell(w)$ of $w$ is the minimal integer $k$ such that $w$ has a word representative of length of~$k$ on the alphabet~$S$. such a word is called a \emph{reduced word representative} of $w$.  \end{definition}

In the sequel, we use the following classical result~\cite{Bou}.
\begin{Prop} Let $(W,S)$ be a \emph{Coxeter system} and $I,J$ be subsets of $S$.  Let $W_I$ and $W_J$ be the subgroups of $W$  generated  by $I$ and $J$ respectively.\\(i) The pairs $(W_I,I)$ and $(W_J,J)$ are Coxeter systems.\\(ii) For every element $w$ which belongs to $W$ there exists a unique element~$\hat{w}$ of minimal length in the double-class $W_JwW_I$. Furthermore there exists $w_1$ in $W_I$ and $w_2$ in $W_J$ such that $w = w_2\hat{w}w_1$ with $\ell(w) = \ell(w_1)+\ell(\hat{w})+\ell(w_2)$.  \label{proppopa} 
\end{Prop}
Note that~$(ii)$ holds when $I$ or $J$ are empty.
%--------------------------------------------------
\subsection{Cross section}
%--------------------------------------------------
Our objective here is to prove Theorem~\ref{thepreserennmonidintro}. We first need to precise the notation used in this theorem. \label{sousectionpresrenmnoid} In all this section, we assume $M$ is a regular irreducible algebraic monoid with a zero element. We denote by~$G$ the unit group of $M$. We fix a maximal torus~$T$ of $G$ and we denote by $\sn$ the Weyl group~$N_G(T)/T$ of~$G$. We denote by $S$ the standard generating set associated with the canonical Coxeter structure of the Weyl group~$W$.

%----------------------------------------
\subsubsection{The Cross Section Lattice}
%----------------------------------------
To describe the generating set of our presentation, we need to introduce the \emph{cross section lattice}, which is related to Green's relations. The latter are classical tools in semigroup theory. Let us recall the definition of Relation~$\mathcal{J}$. The~$\mathcal{J}$-class of an element~$a$ in~$M$ is the double coset~$MaM$. The set~$\mathcal{U}(M)$ of~$\mathcal{J}$-classes carries a natural partial order~$\leq$ defined by~$MaM \leq MbM$ if~$MaM\subseteq MbM$. It turns out that the map~$e\mapsto MeM$ from~$E(M)$ to~$\mathcal{U}(M)$ induces a one-to-one correspondence between the set of~$\sn$-orbits on~$E(\overline{T})$ and the set~$\mathcal{U}(M)$. The existence of this one-to-one correspondence leads to the following definition:
\begin{definition}[cross section lattice] A subset~$\crl$ of~$E(\overline{T})$ is a \emph{cross section lattice} if it is a transversal of~$E(\overline{T})$ for the action of~$\sn$ such that the bijection~$e\mapsto MeM$ from~$\crl$ onto~$\mathcal{U}(M)$ is order preserving.
\end{definition}
It is not immediatly clear that such a cross section lattice exists. Indeed it is, and
\begin{The}\cite[Theorem 9.10]{Put} For every Borel subgroup~$B$ of $G$ that contains $T$, we set $$\crl(B) = \{e\in E(\overline{T})\mid \forall b\in B,\ be = ebe \}.$$ The map~$B\mapsto \crl(B)$ is a bijection between the set of Borel subgroups of~$G$ that contain~$T$ and the set of cross section lattices of~$E(T)$.      
\end{The}   

\begin{Exe}\label{exemplecrossselat} Consider $M = M_n$. Consider the Borel subgroups~$\mathbb{B}$ of invertible upper triangular matrices and~$\mathbb{T}$ the maximal torus of invertible diagonal matrices. Denote by~$e_i$  the diagonal matrix~$\left(\begin{array}{cccccc}Id_i&0\\0&0\end{array}\right)$ of rank~$i$. Then, the set~$\crl(\mathbb{B})$  is~$\{e_0,\ldots, e_n\}$. One has $e_i\leq e_{i+1}$ for every index~$i$.
\end{Exe}
\begin{Rem} \label{remgammainclamb}(i) Let $\Gamma$ be a maximal chain of idempotent elements of $\overline{T}$ and consider the Borel subgroup~$B^+(\Gamma)$ defined in Theorem~\ref{theoemoppborel}. It follows from the definitions that we have~$\Gamma\subseteq \Lambda(B^+(\Gamma))$.\\(ii)\cite[Def.~9.1]{Put} A cross section lattice is a sublattice of $E(\overline{T})$.   
\end{Rem} 

%----------------------------------------
\subsubsection{Normal form and type map}
%----------------------------------------
In order to state the defining relations of our presentation, we turn now to the notion of a {type map}. We fix a Borel subgroup of $G$ that contains~$T$. We write $\crl$ for $\crl(B)$. We consider~$\crl$ as a sublattice of $E(\renn)$ ({\it cf.} Proposition~\ref{propRgenparEetW}).
\begin{Not}\cite{Ren} Let~$e$ belong to~$\crl$.(i) We  set $$\lambda(e) = \{s\in S\mid se = es\}.$$ The map $\lambda : e\mapsto \lambda(e)$ is called the \emph{type map} of the reductive monoid~$M$.\\ (ii) We set $\lambda_\star(e) = \bigcap_{f\geq e}\lambda(f)$ and  $\lambda^\star(e) = \bigcap_{f\geq e}\lambda(f)$.\\ (iii) We set~$\ZZe = \{w\in \sn\mid we = ew \}$, $\ZZse = \{w\in \ZZe \mid we = e\}$ and $\ZZie = \{w\in \ZZe \mid we \neq e\}$. \end{Not}
\begin{Prop}\cite[Lemma 7.15]{Ren} . Then\\ (i)  $\lambda_\star(e) = \{s\in S\mid se = es = e\}$ and $\lambda^\star(e) = \{s\in S\mid se = es \neq e\}$.\\
(ii) The sets~$\ZZe$, $\ZZse$ and $\ZZie$ are the \emph{standard parabolic subgroups} of~$\sn$ generated by the sets $\lambda(e)$ and $\lambda_\star(e)$ and $\lambda^\star(e)$, respectively.  Furthermore, $\ZZe$ is the direct product of $\ZZse$ and $\ZZie$.\label{propopo} \end{Prop}

\begin{Not}\cite{Ren} By Propositions~\ref{proppopa} and~\ref{propopo}, for every~$w$ in~$\sn$ and every~$e,f$ in~$\crl$, each of the sets~$w\ZZe$, $\ZZe w$, $w\ZZse$, $\ZZse w$ and $\ZZe w \ZZf$ has a unique element of minimal length. We denote by $\De,\Ge,\Dse$ and $\Gse$ the set of elements~$w$ of~$\sn$ that are of minimal length in their classes $w\ZZe$, $\ZZe w$, $w\ZZse$ and $\ZZse w$, respectively. Note that the set of elements~$w$ of~$\sn$ that are of minimal length in their double class $\ZZe w\ZZf$ is $\Ge\cap \Df$.
\end{Not}

%----------------------------------------
\subsubsection{Properties of the cross section lattice}
%----------------------------------------
As in previous sections, we fix a Borel subgroup~$B$ of $G$ that contains $T$, and denote by $\Lambda$ the associated cross section lattice contained in~$E(\renn)$. We use the notation~$\mathcal{E}(\Gamma)$ of Section~\ref{sousectionrappelcgt}. We set $\crlo = \crl-\{1\}$. To make the statement of Proposition~\ref{thepreserennmonid} clear we need a preliminary result.

\begin{Lem} Let $e_1,e_2$ lie in $E(T)$ such that $e_1\leq e_2$. There exists $f_1,f_2$ in $\Lambda$ with $f_1\leq f_2$ and $w$ in $W$ such that $wf_1w^{-1} = e_1$ and $wf_2w^{-1} = e_2$. \label{lemtechnorder}
\end{Lem}
\begin{proof}
 Let $\Gamma$ be a maximal chain of $E(\renn)$ that contains $e_1$ and $e_2$. The Borel subgroup $B^+(\Gamma)$ contains the maximal torus~$T$. Therefore, there exists~$w$ in $W$ such that $w^{-1}B^+(\Gamma) w = B$. This implies that $w^{-1}\Lambda(B^+(\Gamma)) w = \crl$. We conclude using Remark~\ref{remgammainclamb}$(i)$. 
\end{proof}
\begin{Lem} \label{lem:cok} Let $h,e$ belong to $\crl$ such that $h \leq e$. Then, $W(h)\cap G(e)\subseteq W_\star(h)$ and $W(h)\cap D(e)\subseteq W_\star(h)$.  
\end{Lem}
\begin{proof} 
Let $w$ lie in~$W(h)\cap G(e)$. We can write $w = w_1w_2 = w_2w_1$ where $w_1$ lies in $W_\star(h)$ and $w_2$ lies in $W^\star(h)$. Since $h\leq e$, we have $\lambda^\star(h)\subseteq\lambda^\star(e)$ and $W^\star(h)\subseteq W^\star(e)$. Since $w$ is assumed to belong to $G(e)$, this implies~$w_2 = 1$. The proof of the second inclusion is similar. 
\end{proof}

\begin{Prop} Let $e,f$ lie in $\crlo$ and $w$ lie in~$\Ge\cap \Df$. There exists $h$ in $\crlo$ with $h\leq e\wedge f$ such that $w$ belongs to $W_\star(h)$ and $ewf = hw = h$.  \label{propconslong}
\end{Prop}
To prove the above proposition, we are going to use the existence of a normal decomposition in~$\renn$: 
\begin{Prop}[\cite{Ren} Section 8.6] \label{fnrenner} For every~$w$ in~$\renn$ there exists a unique triple~$(w_1,e,w_2)$ with~$e\in \crl$,~$w_1\in \Dse$ and~$w_2\in \Ge$ such that~$w = w_1ew_2$. \end{Prop}
Following~\cite{Ren}, we call the triple~$(w_1,e,w_2)$ the \emph{normal decomposition} of $w$.
\begin{proof}[Proof of Proposition~\ref{propconslong}]
Consider the normal decomposition~$(w_1,h,w_2)$ of $ewf$. Then,  $w_1$ belongs to~$\Dsh$ and $w_2$ belongs to~$\Gh$. The element $ewfw^{-1}$ is equal to $w_1hw_2w^{-1}$ and belongs to $E(\renn)$. Since $w_1$ lies in~$\Dsh$, this implies that $w_3 = w_2w^{-1}w_1$ lies in $W_\star(h)$, and that $e\geq w_1hw_1^{-1}$. By Lemma~\ref{lemtechnorder}, there exists $w_4$ in $\sn$ and $e_1,h_1$ in $\crlo$, with $e_1\geq h_1$, such that $w_4e_1w_4^{-1} = e$ and $w_4h_1w_4^{-1} = w_1hw_1^{-1}$. Since $\Lambda$ is a cross section for the action of $\sn$, we have~$e_1 = e$ and~$h_1 = h$. In particular, $w_4$ belongs to~$W(e)$. Since $w_1$ belongs to $\Dsh$, we deduce that there exists $r$ in $W_\star(h)$ such that $w_4 = w_1r$ with $\ell(w_4) = \ell(w_1)+\ell(r)$. Then, $w_1$ lies in $W(e)$. Now, write $w_2 = w''_2w'_2$ where $w''_2$ lies in $W^\star(h)$ and $w'_2$ belongs to $\Gsh$. One has $ewf = w_1w''_2hw'_2$, and $w_1w''_2$ lies in~$\Dh$. By symmetry, we get that $w'_2$ belongs to $W(f)$. Hence, $w^{-1}_1w{w'}^{-1}_2$ is equal to ${w}_3^{-1}w''_2$ and belongs to $W(h)$. By hypothesis $w$ lies in $\Ge\cap\Df$. Then we must have $\ell({w}_3^{-1}w''_2) = \ell(w^{-1}_1)+\ell({w'}^{-1}_2)+\ell(w)$. Since ${w}_3^{-1}w''_2$ belongs to $W(h)$, it follows that $w_1$ and $w'_2$ belong to $W(h)$ too. This implies $w_1 = w'_2 = 1$ and $w = w_3^{-1}w''_2$. Therefore, $ewf = hw''_2 = hw = wh$. Finally,~$w$ belongs to $W_\star(h)$ by Lemma~\ref{lem:cok}.   
\end{proof}

%----------------------------------------
\subsubsection{A presentation for~$\renn$}
%----------------------------------------

\begin{Not}
(i) For each $w$ in $\sn$, we fix a reduced word representative~$\underline{w}$.\\
(ii) We denote by $e\wedge_w f$ the unique letter in~$\crl$ that represents the element $h$ in Proposition~\ref{propconslong}. 
\end{Not}
Note that for~$s$ in~$S$, one has~$\underline{s} = s$.  We recall that $\crl$ is a sub-lattice of $E(\overline{T})$ for the order~$\leq$ defined by $e\leq f$ if $ef = fe = e$.  We are now ready to state a monoid presentation for $\renn$:

\begin{Prop} The Renner monoid has the following monoid presentation whose generating set is~$S\cup \crlo$ and whose defining relations are:\\  
\begin{tabular}{lll}
(COX1)&$s^2 = 1$,&$s\in S$;\\
(COX2)&$|s,t\rangle^m = |t,s\rangle^m$,&$(\{s,t\},m)\in \mathcal{E}(\Gamma)$;\\
(TYM1)&$se = es$,&$e\in \crlo$, $s\in \lambda^\star(e)$;\\
(TYM2)&$se = es = e$,&$e\in \crlo$, $s\in \lambda_\star(e)$;\\
(TYM3)&$e\underline{w}f = e\wedge_wf$,&$e,f\in \crlo$, $w\in G(e)\cap D(f)$. \\
\end{tabular} \label{thepreserennmonid}
\end{Prop}
Note that when $e\leq f$ and $w = 1$, then Relation~(TYM3) becomes $ef = fe = e$. More generally, one has~$e\wedge_1f = e\wedge f$. 
\begin{proof} By remark~\ref{remgenrset} the submonoids~$E(R)$ and $\sn$ generate the monoid $\renn$. As $S$ is a generating set for $\sn$, it follows from the definition of $\crl$ that the set~$S\cup \crlo$ generates $\renn$ as a monoid. Clearly, Relations~(COX1) and~(COX2) hold in~$\sn$, Relations~(TYM1) and~(TYM2) hold in $\renn$. Relations~(TYM3) hold in~$\renn$ by Proposition~\ref{propconslong}. It remains to prove that we obtain a presentation of the monoid~$\renn$. Let~$w$ belong to~$\renn$ with $(w_1,e,w_2)$ as normal form. Consider any word~$\omega$ on the alphabet~$S\cup \crlo$ that represents~$w$. We claim that starting from~$\omega$, one can obtain the word $\underline{w_1}e\underline{w_2}$ using the relations of the above presentation only. This is almost obvious by induction on the number~$j$ of letters of the word~$\omega$ that belong to~$\crlo$.  
The property holds for~$j = 0$ (in this case~$w = w_1$ and~$e = w_2 = 1$) because~(COX1) and~(COX2) are the defining relations of the presentation of~$\sn$. The case $j =1$ is also clear, applying Relations~(COX1),~(COX2),~(TYM1) and~(TYM2). Now, for $j\geq 2$, the case $j$ can be reduced to the case $j-1$ using Relations~(TYM3) (and the other relations).  
\end{proof}

The presentation  in Proposition~\ref{thepreserennmonid} is not minimal; some relations can be removed in order to obtain the presentation stated in Theorem~\ref{thepreserennmonidintro}. Let us introduce a notation used in this theorem:
\begin{Not} If $e$ lies in $\crl$, we denote by $G^\uparrow(e)$ the set $G(e)\cap \left(\bigcap_{f > e}W(f)\right)$. Similarly, we denote by $D^\uparrow(e)$ the set $D(e)\cap \left(\bigcap_{f > e}W(f)\right)$. 
\end{Not}
\begin{Rem}(i) $\left(\cap_{f > e}\lambda(f)\right)\cap\lambda_\star(e) = \emptyset$ by Proposition~\ref{propopo}.\\ (ii) $$G^\uparrow(e) = G(e)\cap W_{\cap_{f > e}\lambda(f)}\textrm{ and }D^\uparrow(e) = D(e)\cap W_{\cap_{f > e}\lambda(f)}.$$
\end{Rem}

The reader may note that for $e\leq f$ one has $\Gue\cap \Duf = \{1\}$. 
\begin{proof}[Proof of Theorem \ref{thepreserennmonidintro}] We need to prove that every relations $e\underline{w}f = e\wedge_wf$ of type~(TYM3) in Proposition~\ref{thepreserennmonid} can be deduced from Relations~(RBI) of Theorem~\ref{thepreserennmonidintro}, using the other common defining relations of type~(COX1), (COX2), (TYM1) and (TYM2). We prove this by induction on the length of $w$. If $\ell(w) = 0$ then $w$ is equal to $1$ and therefore belongs to $\Gue\cap \Duf$. Assume $\ell(w)\geq 1$ and $w$ does not belong to~$\Gue\cap \Duf$. Assume furthermore~$w$ does not lie in~$\Gue$ (the other case is similar). Choose~$e_1$ in $\crlo$ such that $e_1 > e$ and $w$ does not lie in $W(e_1)$. Then, applying Relations~(RIB), we  can transform the word~$e\underline{w}f$ into the word~$ee_1\underline{w}f$. Using relations~(COX2), we can transform the word~$\underline{w}$ into a word~$\underline{w_1}\,\underline{w_2}$ where $w_1$ belongs to $W(e_1)$ and $w_2$ belongs to $\Gueun$. Then, applying Relations~(COX2) and~(TYM1), we can transform the word~$ee_1\underline{w}f$ into the word~$e\underline{w_1}e_1\underline{w_2}f$. By hypothesis on~$w$, we have~$w_2\neq 1$ and, therefore, $\ell(w_1)< \ell(w)$. Assume $w_2$ belongs  to $\Duf$. We can apply Relation~(RIB) to transform~$e\underline{w_1}e_1\underline{w_2}f$ into $e\underline{w_1}(e_1\wedge_{w_2}f)$. Using relations (COX2), we can transform $\underline{w_1}$ into a word $\underline{w'_1}\,\underline{w''_1}\,\underline{w'''_1}$ with $w'''_1$ in $W_\star(e_1\wedge_{w_2}f)$, $w''_2$ in $W^\star(e_1\wedge_{w_2}f)$ and $w'_1$ in $D(e_1\wedge_{w_2}f)$. Then~$e\underline{w_1}(e_1\wedge_{w_2}f)$ can be transformed into~$e\underline{w'_1}(e_1\wedge_{w_2}f)\underline{w''_1}$.  Since $\ell(w'_1)\leq\ell(w_1)<\ell(w)$, we can apply an induction argument to transform the word~$e\underline{w_1}(e_1\wedge_{w_2}f)$ into the word $e\wedge_{w_1}(e_1\wedge_{w_2}f)\underline{w''_1}$. Now, by the unicity of the normal decomposition, $w''_1$ as to belong to $W_\star(e\wedge_{w_1}(e_1\wedge_{w_2}f))$. Therefore we can transform~$e\wedge_{w_1}(e_1\wedge_{w_2}f)\underline{w''_1}$ into $e\wedge_{w_1}(e_1\wedge_{w_2}f)$ using Relations~(TYM2).   Note that the letters~$e\wedge_{w_1}(e_1\wedge_{w_2}f)$ and $e\wedge_{w}f$ has to be equal as they represent the same element in~$\crl$. Assume finally that $w_2$ does not belong to $\Duf$. By similar arguments we can, applying Relations~(COX2) and~(TYM1), transform the word~$e\underline{w_1}e_1\underline{w_2}f$ into a word~$e\underline{w_1}e_1\underline{w_{3}}f_1\underline{w_{4}}$ where $f_1 > f$ in~$\crlo$ and~$w_2 = w_{3}w_{4}$ with $w_{3}$ in $\Dufun$ and  $w_{4}$ in $W(f_1)$. At this stage we are in position to apply Relation~(RBI). Thus, we can transform the word~$e\underline{w_1}e_1\underline{w_2}f$ into~$e\underline{w_1}\,(e_1\wedge_{w_3}f_1)\,\underline{w_{4}}f$. Since we have~$\ell(w_1)+\ell(w_4) < \ell(w)$ we can proceed as in the first case to conclude.
\end{proof}

%--------------------------------------------------
\section{Some particular Renner monoids}
%--------------------------------------------------
Here we focus on some special Renner monoids considered in~\cite{LiRe,Li1,Li2}. In each case, we provide an explicit monoid presentation using the general presentation obtained in Section~1. 

%--------------------------------------------------
\subsection{The rook monoid}
%--------------------------------------------------
Consider $M = M_n$ and choose $\mathbb{B}$ for Borel subgroup (see Example~\ref{exemplecrossselat}). In this case, the Weyl group is the symmetric group $S_n$. Its generating set $S$ is $\{s_1,\cdots, s_{n-1}\}$ where $s_i$ is the transposition matrix corresponding to~$(i,i+1)$. The cross section section lattice~$\crl = \{e_0,\cdots,e_{n-1},e_n\}$ is linear (we have $e_j\leq e_{j+1}$ for every $j$). For every~$j$ one has $\lambda_\star(e_j) =\{s_i\mid j+1\leq i\}$ and $\lambda^\star(e_j) = \{s_i\mid i\leq j-1\}$. In particular,~$\Guei\cap \Duei = \{1,s_i\}$, and for $i\neq j$ we have~$\Guei\cap \Duej = \{1\}$. 

\begin{figure}[ht]
\begin{picture}(250,75)
\put(18,0){\includegraphics[scale = 0.6]{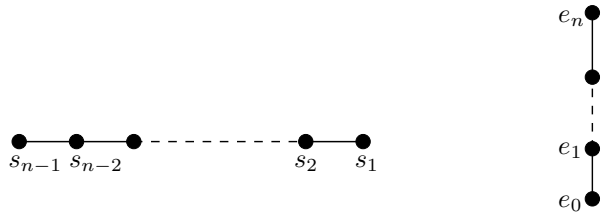}}
\put(148,15){$s_1$} \put(125,15){$s_2$}\put(40,15){$s_{n-2}$}\put(17,15){$s_{n-1}$}
\put(225,0){$e_0$} \put(225,20){$e_1$} \put(225,70){$e_n$}
\end{picture}
\caption{Coxeter graph and Hasse diagram for $M_n$.}\label{fig:hassecrlAn}
\end{figure}

Therefore, we recover the monoid presentation of the rook monoid~$\renn$ stated in~\cite{God}: the generating set is~$\{s_1,\ldots,s_{n-1},e_0,\ldots,e_{n-1}\}$ and the defining relations are $$
\begin{array}{rcll} s_i^2&=&1, &1\leq i \leq n-1;\\s_is_j&=&s_js_i,&1\leq i,j \leq n-1\textrm{ and }|i-j|\geq 2;\\s_is_{i+1}s_i&=&s_{i+1}s_is_{i+1}, &1\leq i \leq n-1;\\e_js_i&=&s_ie_j&1\leq i<j \leq n-1;\\e_js_i&=&s_ie_j = e_j&0\leq j<i\leq n-1;\\e_ie_j&=&e_je_i = e_{\min(i,j)}&0\leq i,j \leq n-1;\\e_is_ie_i&=&e_{i-1} &1\leq i \leq n-1.\end{array}$$ 

%--------------------------------------------------
\subsection{The Sympletic Algebraic Monoid}
%--------------------------------------------------
Let $n$ be\label{simplectmonoid} a positive even integer and $Sp_n$ be the \emph{Symplectic Algebraic Group} \cite[page 52]{Hum}: assume $\ell$ lies in $\mathbb{N}$, and consider the matrix~$J_{\ell} = \left(\begin{array}{ll} &1\\1\end{array}\right)$ in $M_{\ell}$. Let $J = \left(\begin{array}{cc}0&J_{\ell}\\-J_{\ell}&0\end{array}\right)$ in $M_n$, where $n = 2\ell$. Then,~$Sp_n$ is equal to~$\{A\in M_n\mid A^tJA = J\}$, where $A^t$ is the transpose matrix of $A$. We set $M = \overline{K^{\scriptscriptstyle \times}Sp_n}$. This monoid is a regular monoid with $0$ whose associated reductive algebraic unit group is~$K^{\scriptscriptstyle \times}Sp_n$. It is called the \emph{Symplectic Algebraic Monoid} \cite{LiRe}. Let~$\mathbb{B}$ be the Borel subgroup of $GL_n$ as defined in Example~\ref{exemplecrossselat}, and set $B = K^{\scriptscriptstyle \times}(\mathbb{B}\cap Sp_n)$. This is a Borel subgroup of the unit group of $M$. It is shown in~\cite{LiRe}  that the cross section lattice~$\crl$ of $M$ is $\{e_0,e_1,\cdots ,e_{\ell}, e_n\}$ where the elements $e_i$ correspond to the matrices of $M_n$ defined in Example~\ref{exemplecrossselat}. In particular the cross section lattice~$\crl$ is linear. In this case, the Weyl group is a Coxeter group of type~$B_\ell$. In other words, the group~$W$ is isomorphic to the subgroup of $S_{n}$ generated by the permutation matrices $s_1,\cdots,s_\ell$ corresponding to $(1,2)(n-1,n)$, $(2,3)(n-2,n-1)$, $\cdots$, $(\ell-1,\ell) (\ell+1,\ell+2)$, and~$(\ell,\ell+1)$, respectively.  One has $\lambda_\star(e_i) = \{s_{i+1},\cdots, s_{\ell}\}$  and $\lambda^\star(e_i) = \{s_1,\cdots, s_{i-1}\}$. Therefore, $\Guel\cap \Duel = \{1,s_\ell, s_\ell s_{\ell-1}s_\ell\}$, and for $i$ in~$\{1,\cdots, \ell-1\}$ one has $\Guei\cap \Duei = \{1,s_i\}$. A direct calculation proves that $e_is_ie_i = s_i e_{i-1}$ for every $i$, and $e_\ell s_\ell s_{\ell-1}s_\ell e_\ell =  e_{\ell-2}$. 
\begin{figure}[ht]
\begin{picture}(250,75)
\put(18,0){\includegraphics[scale = 0.6]{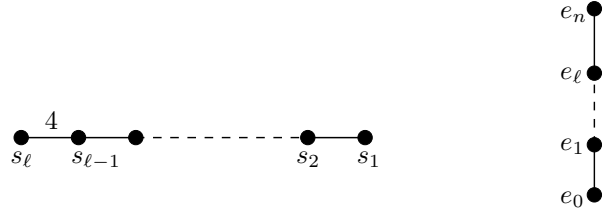}}
\put(148,15){$s_1$} \put(125,15){$s_2$}\put(40,15){$s_{\ell-1}$}\put(30,28){$4$}\put(17,15){$s_{\ell}$}
\put(225,0){$e_0$} \put(225,20){$e_1$} \put(225,47){$e_\ell$}\put(225,70){$e_n$}
\end{picture}
\caption{Coxeter graph and Hasse diagram for $Sp_n$.}\label{fig:hassecrlBn}
\end{figure}
Hence, a monoid presentation of~$\renn$ is given by the generating set~$\{s_1,\ldots,s_{\ell},e_0,\ldots,e_{\ell}\}$ and the defining relations$$ \begin{array}{rcll} s_i^2&=&1, &1\leq i \leq \ell;\\s_is_j&=&s_js_i,&1\leq i,j \leq \ell\textrm{ and }|i-j|\geq 2;\\s_is_{i+1}s_i&=&s_{i+1}s_is_{i+1}, &1\leq i \leq \ell-2;\\s_\ell s_{\ell-1}s_\ell s_{\ell-1}&=&s_{\ell-1}s_\ell s_{\ell-1} s_\ell;\\e_js_i&=&s_ie_j&1\leq i<j \leq \ell;\\e_js_i&=&s_ie_j = e_j&0\leq j<i\leq \ell;\\e_ie_j&=&e_je_i = e_{\min(i,j)}&0\leq i,j \leq \ell;\\e_is_ie_i&=&e_{i-1} &1\leq i \leq  \ell;\\e_\ell \, s_{\ell}s_{\ell-1} s_\ell\, e_\ell&=&e_{\ell-2} .\end{array}$$ 
%--------------------------------------------------
\subsection{The Special Orthogonal Algebraic Monoid}
%--------------------------------------------------
Let $n$ be a positive integer and $J_n$ be defined as in Section~\ref{simplectmonoid}. The \emph{Special Orthogonal Group}~${\bf SO}_n$ is defined as ${\bf SO}_n  = \{A\in SL_n\mid g^T J_n g = J_n\}$. The group $K^{\scriptscriptstyle \times}\, {\bf SO}_n$ is a connected reductive group. Following~\cite{Li1,Li2}, we define the \emph{Special Orthogonal Algebraic Monoid} to be the Zariski closure~$M = \overline{K^{\scriptscriptstyle \times}{\bf SO}_n}$ of $K^{\scriptscriptstyle \times}{\bf SO}_n$. This is an algebraic monoid~\cite{Li1,Li2}, and $B = \mathbb{B} \cap M$ is a Borel subgroup of its unit group.   In this case, the cross section lattice depends on the parity of $n$. Furthermore, the Weyl group is a Coxeter group whose type depends on the parity of $n$ too.

\begin{figure}[ht]
\begin{picture}(250,100)
\put(18,0){\includegraphics[scale = 0.6]{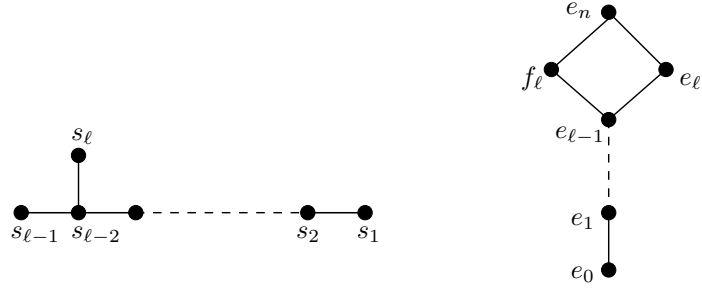}}
\put(148,15){$s_1$} \put(125,15){$s_2$}\put(40,15){$s_{\ell-2}$}\put(17,15){$s_{\ell-1}$}\put(40,52){$s_\ell$}
\put(229,0){$e_0$} \put(229,20){$e_1$} \put(223,53){$e_{\ell-1}$} \put(210,73){$f_{\ell}$} \put(270,73){$e_\ell$} \put(227,100){$e_n$}
\end{picture}
\caption{Coxeter graph and Hasse diagram for ${\bf SO}_{2\ell}$.}\label{fig:hassecrlDn}
\end{figure}
Assume $n = 2\ell$ is even. In this case, $W$ is a Coxeter group of type~$D_\ell$. The standard generating set of $W$ is $\{s_1,\cdots, s_{\ell}\}$ where for $1\leq i\leq \ell-1$, the element~$s_i$ is the permutation matrix associated with~$(i,i+1)(n-i,n-i+1)$, and $s_\ell$ is the permutation matrix associated with~$(\ell-1,\ell+1)(\ell,\ell+2)$. It is shown in~\cite{Li2} that the cross section~$\crl$ is equal to~$\{e_0,e_1,\cdots,e_{\ell},f_\ell,e_n\}$.. The elements $e_i$ correspond to the matrices of $M_n$ defined in Example~\ref{exemplecrossselat}; the element~$f_{\ell}$ is the diagonal matrix~$e_{\ell+1}+e_{\ell-1}-e_\ell$. The Hasse diagram of $\crl$ is as represented in Figure~\ref{fig:hassecrlDn}.  For~$j$ in $\{0,\ldots,\ell-2\}$ one has $\lambda_\star(e_j) =\{s_i\mid j+1\leq i\}$ and $\lambda^\star(e_j) = \{s_i\mid i\leq j-1\}$. Furthermore, one can verified that $$\lambda_\star(e_{\ell-1}) = \lambda_\star(f_{\ell})  = \lambda_\star(e_{\ell})= \emptyset,$$ $$\lambda^\star(e_{\ell-1}) = \lambda^\star(f_{\ell})=  \{s_i\mid i\leq \ell-2\},$$ $$\lambda^\star(e_{\ell})=  \{s_i\mid i\leq \ell-1\},$$ $$\lambda^\star(f_{\ell})= \{s_i\mid i\neq\ell-1\}.$$ 
Therefore, for $i$ in~$\{1,\cdots, \ell-2\}$  one has $\Guei\cap \Duei = \{1,s_i\}$. Furthermore, 
$$\Guell\cap \Duell = \{1\} \textrm{ and }$$ $$\Gufl\cap \Dufl = \{1,s_{\ell-1}\}\ \ \ ;\ \ \ \Guel\cap \Duel = \{1,s_\ell\};$$ $$\Guel\cap \Dufl = \{1, s_\ell s_{\ell-2}s_{\ell-1}\}\ \ \ ;\ \ \ \Gufl\cap \Duel = \{1, s_{\ell-1} s_{\ell-2} s_\ell\}.$$  The monoid~$\renn$ has a presentation with~$\{s_1,\ldots,s_{\ell},e_0,\ldots,e_{\ell},f_\ell\}$ for generating set and  
$$
\begin{array}{rcll} s_i^2&=&1, &1\leq i \leq \ell;\\s_is_j&=&s_js_i,&1\leq i,j \leq \ell\textrm{ and }|i-j|\geq 2;\\s_is_{i+1}s_i&=&s_{i+1}s_is_{i+1}, &1\leq i \leq \ell-2;\\s_\ell s_{\ell-2}s_\ell &=&s_{\ell-2}s_\ell s_{\ell-2};\\e_js_i&=&s_ie_j,&1\leq i<j \leq \ell;\\e_js_i&=&s_ie_j = e_j,&0\leq j<i\leq \ell;\\e_ie_j&=&e_je_i = e_{\min(i,j)},&0\leq i,j \leq \ell;\\f_\ell e_\ell&=&e_\ell f_\ell = e_{\ell-1};\\e_is_ie_i&=&e_{i-1}, &1\leq i \leq  \ell-1;\\e_\ell s_\ell e_\ell &=&f_\ell s_{\ell-1} f_\ell = e_{\ell-2};\\
e_\ell \, s_{\ell}s_{\ell-2} s_{\ell-1}\, f_\ell&=&f_\ell \, s_{\ell-1}s_{\ell-2} s_{\ell}\, e_\ell = e_{\ell-3} .\end{array}$$ 
for defining relations.
Assume $n = 2\ell+1$ is odd.  In that case, $W$ is a Coxeter group of type~$B_\ell$. It is shown in~\cite{Li1} that the cross section lattice is linear as in the case of the Symplectic Algebraic Monoid. It turns out that the Renner monoid of~${\bf SO}_{2\ell+1}$ is isomorphic to the Renner monoid of Symplectic Algebraic Monoid~$\overline{K^{\scriptscriptstyle \times} Sp_{2\ell}}$, and that we obtain the same presentation than in the latter case.
%--------------------------------------------------
\subsection{More examples: Adjoint Representations}
%--------------------------------------------------
Let~$G$ be a simple algebraic group, and denote by~$\mathfrak{g}$ its Lie algebra. Let $M$ be the algebraic monoid~$\overline{K^{\scriptscriptstyle \times} Ad(G)}$ in $End(\mathfrak{g})$. The cross section lattice of~$M$ and the type map of~$M$ has been calculated for each Dynkin diagram (see~\cite[Sec.~7.4]{Ren}).  Therefore one can deduce a monoid presentation for each of the associated Renner monoid. 
%--------------------------------------------------
\subsection{More examples: $\mathcal{J}$-irreducible algebraic monoids}
%--------------------------------------------------
In~\cite{RePu}, Renner and Putcha consider among regular irreducible algebraic monoids those who are \emph{$\mathcal{J}$-irreducible}, that is those whose cross section lattices have a unique minimal non-zero element. It is easy to see that the $\mathcal{J}$-irreducibility property  is related to the existence of irreducible rational representations \cite[Prop.~4.2]{RePu}. Renner and Putcha determined the cross section lattice of those~$\mathcal{J}$-irreducible that arise from special kind of dominant weigths~\cite[Fig.~2,3]{RePu}. Using~\cite[Theorem~4.13]{RePu}, one can deduce the associated type maps and therefore a monoid presentation of each corresponding Renner monoids.   

%--------------------------------------------------
%--------------------------------------------------
\section{A length function on $\renn$}
%--------------------------------------------------
%--------------------------------------------------

In this section we extend  the length function defined in~\cite{God} to any Renner monoid. In all this section, we assume $M$ is a regular irreducible algebraic monoid with a zero element. We denote by~$G$ the unit group of $M$. We fix a maximal torus~$T$ of $G$ and we denote by $\sn$ the Weyl group~$N_G(T)/T$ of~$G$. We denote by $S$ the standard generating set associated with the canonical Coxeter structure of the Weyl group~$W$. We fix a Borel subgroup~$B$ of $G$ that contains $T$, and we denote by~$\crl$ the associated cross section lattice contained in~$\renn$. As before we set $\crlo = \crl-\{1\}$.

%--------------------------------------------------
\subsection{Minimal word representatives}
%--------------------------------------------------
 The definition of the length function on $W$ and of a reduced word has been recalled in Section~\ref{sousectionrappelcgt}. If $X$ is a set, we denote by $X^*$ the set of finite words on $X$.

\begin{definition}\label{deflenfun} (i) We set~$\ell(s) = 1$ for~$s$ in~$S$ and~$\ell(e) = 0$ for~$e$ in~$\crl$. Let~$x_1,\ldots, x_k$ be in~$S\cup\crlo$  and consider the word~$\omega  =  x_1\cdots x_k$. Then, the \emph{length} of the word~$\omega $ is the integer~$\ell(\omega )$ defined by~$\ell(\omega ) = \sum_{i = 1}^k\ell(x_i)$.\\
(ii) The \emph{length} of an element~$w$ which belongs to $\renn$ is the integer~$\ell(w)$ defined by $$\ell(w) =  \min\left\{\ell(\omega), \omega \in (S\cup\crlo)^*\mid \omega \textrm{ is a word representative of }w \right\}.$$ \end{definition}
The following properties are direct consequences of the definition. 
\begin{Prop} \label{propproplenght}Let~$w$ belong to~$\renn$.\\ (i) The length function~$\ell$ on~$\renn$ extends the length function~$\ell$ defined on~$\sn$.\\(ii) If~$\ell(w) = 0$ then~$w$ lies in~$\crl$.\\(iii) If~$s$ lies in~$S$ then~$|\ell(sw)-\ell(w)|\leq 1$.\\(iv) If~$w'$ belongs to~$\renn$, then~$\ell(ww')\leq \ell(w)+\ell(w')$.\end{Prop}

\begin{proof} Point~$(i)$ and~$(ii)$ are clear (the letters of every representative word of an element in~$\sn$ are in~$S$). If~$w$ lie in~$\renn$ and~$s$ lie to~$S$, then~$\ell (sw)\leq \ell(w)+1$. Since~$w = s^2w = s (sw)$, the inequality~$\ell(w)\leq \ell(sw)+1$ holds too. Point~$(iii)$ follows, and Point~$(iv)$ is a direct consequence of~$(iii)$.
\end{proof}

\begin{Prop}\label{propproplenght2}
Let~$w$ belong to~$\renn$.\\(i) If $(w_1,e,w_2)$ is the normal decomposition of $w$, then $\ell(w) = \ell(w_1)+\ell(w_2)$.\\(ii) If $\omega_1,\omega_2$ are two representative words of $w$ on $S\cup \crlo$ such that the equalities~$\ell(w) = \ell(\omega_1) = \ell(\omega_2)$ hold, then using the defining relations of the presentation of $\renn$ in Theorem~\ref{thepreserennmonid}, we can transform $\omega_1$ into $\omega_2$ without increasing the length. 
\end{Prop}
\begin{proof}
(i) Let $\omega$ be a representative word $w$ on $S\cup \crlo$ such that $\ell(w) = \ell(\omega)$. It is clear that we can repeat the argument of the proof of Theorem~\ref{thepreserennmonid} without using the relation $(BR1)$. Therefore~$\ell(\omega)\geq \ell(w_1ew_2) = \ell(w_1)+\ell(w_2)\geq(w)$.\\
(ii) is a direct consequence of the proof of $(i)$. \end{proof}

\begin{Cor} Let $w$ lie in~$\renn$ and $e$ belongs to~$\crlo$. Denote by $(w_1,f,w_2)$ the normal decomposition of $w$.\\ (i) One has~$\ell(we)\leq \ell(w)$ and~$\ell(ew)\leq \ell(w)$.\\ (ii) $\ell(we) = \ell(w)$ if and only if the normal decomposition of $we$ is $(w_1,e\land f, w_2)$. Furthermore, in this case, $w_2$ lies in $W^\star(e)$. 
\end{Cor}
\begin{proof} (i) is a direct consequence of the definition of the length and of Proposition~\ref{propproplenght2}(i) : $\ell(we) = \ell(w_1fw_2e)\leq \ell(w_1)+0 +\ell(w_2)+0 = \ell(w)$. The same arguments prove that~$\ell(ew)\leq \ell(w)$.\\(ii) Decompose~$w_2$ as a product~$w'_2w''_2w'''_2$ where $w'''_2$ lies in $W_\star(f)$, $w''_2$ lies in $W^\star(f)$, $w'_2$ lies in $D(f)$ and $\ell(w_2) = \ell(w'_2)+\ell(w''_2)+ \ell(w'''_2)$. Then, $we = w_1fw_2e = w_1fw'_2ew''_2 = w_1(f\land_{w'_2}e)w''_2$. In particular, $\ell(we)\leq \ell(w_1)+\ell(w''_2)$. Assume $\ell(we) = \ell(w)$. We must have $w'_2 = w'''_2 = 1$. The element~$w''_2$, that is $w_2$, must belong to $D^\star(f\land_1 e) = D^\star(f\land e)$, and the element~$w_1$ must belong to $G^\star(f\land e)$. In particular, $w_2$ lies in $W^\star(e)$. Furthermore, $w_2$ lies in $D(f\land e)$ since $\lambda^\star(f\land e)\subseteq \lambda(f)$ by Proposition~\ref{propopo}(i).  Conversely, if the the normal decomposition of $we$ is $(w_1,e\land f, w_2)$, then  $\ell(we) = \ell(w_1)+\ell(w_2) = \ell(w)$. 
 
\end{proof}

%-------------------------------------------------
\subsection{Geometrical formula}
%-------------------------------------------------
In Proposition~\ref{propgeominterpr} below we provide a geometrical formula for the length function~$\ell$ defined in the previous section. This formula extends naturally the geometrical definition of the length function on a Coxeter group. Another length function on Renner monoids has already been defined and investigated~\cite{Sol5,PePuRe,Ren}. This length function has nice properties, which are similar to the ones in Propositions~\ref{propproplenght}, \ref{propproplenght2}(i) and~\ref{propgeominterpr}. This alternative length function has been firstly  introduced by Solomon~\cite{Sol5} in the special case of rook monoids in order to verify a combinatorial formula that generalizes Rodrigues formula~\cite{Rod}. That is why we call this length function the \emph{Solomon length function} in the sequel. We proved in~\cite{God} that our length function for the rook monoid verifies the same combinatorial formula. We also proved in~\cite{God} that in the case of the rook monoid, our presentation of~$\renn$ and our length function are related to the Hecke algebra. We believe this is still true in the general case.  
\begin{Lem} \label{lemtechimpo} Let~$w$ belong to~$\renn$ and denote by~$(w_1,e,w_2)$ its normal decomposition. Let~$s$ be in~$S$.\\(i) We have one of the two following cases:\\\indent (a) there exists~$t$ in~$\lambda_\star (e)$ such that~$sw_1 = w_1t$. In this case, one has~$sw = w$;\\\indent (b) the element~$sw_1~$ lies in~$\Dse$ and~$(sw_1,e,w_2)$ is the normal decomposition of~$sw$.\\(ii) Denote by~$\tilde{l}$ the Solomon length function on~$\renn$ . Then, $$\ell(sw) - \ell(w) = \tilde{l}(sw) -\tilde{l}(w).$$     
\end{Lem}  
\begin{proof}~$(i)$ If~$sw_1$ lies in~$\Dse$, then by Theorem~\ref{fnrenner}, the triple~$(sw_1,e,w_2)$ is the normal decomposition of~$sw$. Assume now that~$sw_1$ does not belong to~$\Dse$. In that case,~$e$ cannot be equal to~$1$. Since~$w_1$ belongs to~$\Dse$, by the exchange lemma, there exists~$t$ in~$\lambda_(e)$  such that~$sw_1 = w_1t$. Therefore,~$sw = sw_1ew_2 = w_1tew_2 = w_1ew_2= w$.\\
$(ii)$ The Solomon length~$\tilde{l}(w)$ of an element~$w$ in~$\renn$ can be defined by the formula~$\tilde{l}(w) = \ell(w_1)-\ell(w_2)+\tilde{\ell}_e$ where~$(w_1,e,w_2)$ is the normal decomposition of~$w$ and~$\tilde{\ell}_e$ is a constant that depends on~$e$ only \cite[ Definition~4.1]{PePuRe}. Therefore the result is a direct consequence of~$(i)$. 
\end{proof}   

As a direct consequence of Lemma~\ref{lemtechimpo}$(ii)$ and~\cite[Theorem~8.18]{Ren} we get Proposition~\ref{propintro}.

\begin{Prop}\label{propgeominterpr} Let~$w$ belong to~$\renn$, and~$(w_1,e,w_2)$ be its normal decomposition. $$\ell(w) = \dim(B w_1e B) - \dim(B ew_2 B).$$\end{Prop}
When~$w$ lies in~$S_n$, one has~$e = w_2 = 1$, and we recover the well-known formula~$\ell(w)  = \dim(BwB) - \dim(B)$.
\begin{proof} By~\cite[Section~4]{PePuRe}, for every normal decomposition~$(v_1,f,v_2)$ we have the equality~$\dim(B v_1fv_2 B) = \ell(v_1)-\ell(v_2)+k_f$, where~$k_f$ is a constant that depends on~$f$ only. Therefore,$$\dim(B w_1e B) - \dim(B ew_2 B) = \ell(w_1)+k_e - (-\ell(w_2)+k_e) = \ell(w).$$
\end{proof}

%-------------------------- 
\

\end{document}